\documentclass[preprint,11pt]{article}
\usepackage{amsmath,amssymb}
\usepackage{multirow}
\usepackage{enumerate}
\usepackage{array}
\usepackage{enumerate}
\usepackage{algorithm}
\usepackage{algpseudocode}
\usepackage{algorithmicx}
\usepackage{graphicx}

\newtheorem{theo}{Theorem}[section]

\newtheorem{coro}[theo]{Corollary}
\newtheorem{lemm}[theo]{Lemma}
\newtheorem{prop}[theo]{Proposition}

\begin{document}

\begin{center}
{\Large \bf {Algebraic computation of genetic patterns related to three-dimensional evolution algebras\footnote{Preliminary version. Accepted manuscript for publication in {\em Applied Mathematics and Computation}. DOI: 10.1016/j.amc.2017.05.045.}}}
\end{center}

\begin{center}
{\large \em O. J. Falc\'on$^1$,\ R. M. Falc\'on$^2$\, J. N\'u\~nez$^{3}$}
\vspace{0.5cm}

{\small $^{1, 3}$ Faculty of Mathematics, University of Seville, Spain.\\
$^2$School of Building Engineering, University of Seville, Spain.\\
E-mail: {\em $^1$oscfalgan@yahoo.es,\ $^2$rafalgan@us.es,\ $^3$jnvaldes@us.es}}
\end{center}

\vspace{0.5cm}

\noindent {\large \bf Abstract.} The mitosis process of an eukaryotic cell can be represented by the structure constants of an evolution algebra. Any isotopism of the latter corresponds to a mutation of genotypes of the former. This paper uses Computational Algebraic Geometry to determine the distribution of three-dimensional evolution algebras over any field into isotopism classes and hence, to describe the spectrum of genetic patterns of three distinct genotypes during a mitosis process. Their distribution into isomorphism classes is also determined in case of dealing with algebras having a one-dimensional annihilator.

\vspace{0.5cm}


\section{Introduction}

Evolution algebras were introduced by Tian and Vojtechovsky \cite{Tian2008, Tian2006} to simulate algebraically the self-reproduction of alleles in non-Mendelian Genetics. In the last years, these algebras have been widely studied without probabilistic restrictions on their structure constants \cite{Cabrera2016, Camacho2013a, Casas2011, Casas2013, Dzhumadildaev2016, Labra2014, Ladra2013, Ladra2014, Omirov2015}. A main problem to be solved in this regard is the distribution of evolution algebras into isomorphism classes. Besides, Holgate and Campos \cite{Campos1987, Holgate1966} exposed the importance in Genetics of considering also the distribution of such algebras under isotopism classes, because they constitute a way to formulate algebraically the mutation of genotypes in the inheritance process.

\vspace{0.1cm}

The mentioned distribution into isomorphism classes has already been dealt with for two-dimensional evolution algebras over the complex field \cite{Camacho2013, Casas2014} and for nilpotent evolution algebras of dimension up to four over arbitrary fields \cite{Hegazi2015}. More recently, the authors \cite{Falcon2016_Evol} have characterized the isomorphism classes of two-dimensional evolution algebras over arbitrary fields and have established their distribution into isotopism classes. The latter gives rise to the spectrum of genetic patterns of two distinct genotypes during mitosis of eukaryotic cells. More specifically, if $u$ and $v$ are two distinct non-zero elements of the algebra under consideration, then the mentioned spectrum is formed by the next four genetic patterns: $(0,0)$, for which no offspring exists; $(u,0)$, for which exactly one of the two genotypes does not produce offspring; $(u,u)$, for which the offspring has always the same genotype, whatever the initial one is; and $(u,v)$, for which the genotype of the offspring depends directly on that of the cell parent. This spectrum was obtained in \cite{Falcon2016_Evol} by means of distinct aspects on Computational Algebraic Geometry, all of them based on the fact that the algebraic law defined by the structure constants of any evolution algebra, together with the relations among basis vectors described by any isotopism of algebras, constitutes the set of generators of an ideal of polynomials whose reduced Gr\"obner basis establishes the algebraic relations that must hold, up to mutation, the genetic patterns of the mitosis process. This paper also focuses on the computation of such relations, particularly in case of dealing with three-dimensional evolution algebras. This enables us to distribute these algebras into isotopism classes, whatever the base field is, and describe mathematically the spectrum of genetic patterns of three distinct genotypes during a mitosis process. The distribution of such algebras into isomorphism classes is also determined in case of dealing with algebras having a one-dimensional annihilator.

\section{Preliminaries}

In this section we expose some basic concepts and results on evolution algebras, isotopisms and Computational Algebraic Geometry that are used throughout the paper. For more details about these topics we refer the reader to the respective manuscripts of Tian \cite{Tian2008}, Albert \cite{Albert1942} and Cox et al. \cite{Cox1998}.

\subsection{Evolution algebras}

A {\em gene} is the molecular unit of hereditary information. This consists of {\em deoxyribonucleic acid} ({\em DNA}), which contains the code to determine the attributes or {\em phenotypes} that characterize each organism. Genes are disposed in sequential order giving rise to long strands of DNA called {\em chromosomes}. Genes related to a given phenotype can have distinct forms, which are called {\em alleles}. They always appear at the same position in chromosomes and constitute the {\em genotype} of the organism with respect to such a phenotype.

\vspace{0.1cm}

In {\em eukaryotic} cells, DNA is primordially contained in the {\em nucleus}, although it also appears in organelles as mitochondria and chloroplasts, which are located in the cytoplasm. {\em Mitosis} is an asexual form of reproduction that consists of the division of an eukaryotic cell into two daughter cells in such a way that the nuclear genetic material of the former duplicates giving rise, up to rare mutation, to homologous nuclear chromosomes. In the final stage, the two daughter cells split having each resulting cell its corresponding copy of nuclear genetic material, whereas the extra-nuclear genetic material in the cytoplasm of the parent cell is randomly distributed between both of them by {\em vegetative division}. There exist distinct probabilistic laws that regulate the theoretical influence of all of this genetic material in the genotype of the offspring. Tian and Vojtechovsky \cite{Tian2006} introduced evolution algebras to represent mathematically these laws. Specifically, an {\em evolution algebra} defined on a set $\beta=\{e_1,\ldots, e_n\}$ of distinct genotypes with respect to a given phenotype of an asexual organism is an $n$-dimensional algebra of basis $\beta$ over a field $\mathbb{K}$ such that $e_ie_j=0$, if $i\neq j$, and $e_ie_i=\sum_{j=1}^n c_{ij}e_j$, for some $c_{i1},\ldots,c_{in}\in \mathbb{K}$. The elements $c_{ij}$ are called the {\em structure constants} of the algebra. Here, the product $e_ie_j=0$, for $i\neq j$, is due to the uniqueness of the genotype of the parent cell; the product $e_ie_i$ represents self-replication in the mitosis process; and each structure constant $c_{ij}$ constitutes the probability, due to vegetative division, that the genotype $e_i$ becomes the genotype $e_j$ in the next generation.

\subsection{Isotopisms of evolution algebras}

Two $n$-dimensional algebras $A$ and $A'$ are {\em isotopic} \cite{Albert1942} if there exist three non-singular linear transformations $f$, $g$ and $h$ from $A$ to $A'$ such that $f(u)g(v)=h(uv)$, for all $u,v\in A$. The triple $(f,g,h)$ is called an {\em isotopism} between $A$ and $A'$. If $f=g$, then the triple $(f,f,h)$ is called a {\em strong isotopism} and the algebras are said to be {\em strongly isotopic}. If $f=g=h$, then the isotopism constitutes an {\em isomorphism}, which is denoted by $f$ instead of $(f,f,f)$. To be isotopic, strongly isotopic or isomorphic are equivalence relations among algebras. Throughout the paper, we refer the former and the latter as $\sim$ and $\cong$, respectively.

\vspace{0.1cm}

Isotopisms of evolution algebras can be interpreted as mutations of the genetic material of parent and daughter cells in the mitosis process with respect to a given phenotype. Specifically, if $(f,g,h)$ is an isotopism between two $n$-dimensional evolution algebras $A$ and $A'$, then $f$ and $g$ represent the respective possible mutation of each one of the two homologous chromosomes in which the nuclear genetic material of the parent cell duplicates during the first part of the mitosis process, whereas $h$ represents a possible mutation of the genotype of the offspring in the final step of the process. For each $\alpha\in\{f,g,h\}$, the corresponding expression $\alpha(e_i)=\sum_{j=1}^n a_{ij}e_j$ involves the genotype $e_i$ to mutate to $e_j$ with probability $a_{ij}$. Since $A'$ is also an evolution algebra, the mitosis process only finishes if the genotypes of both homologous chromosomes that have been created after the mutations $f$ and $g$ coincide. If these genotypes do not coincide, then there is no offspring.

\vspace{0.1cm}

Hereafter, $\mathcal{E}_n(\mathbb{K})$ denotes the set of $n$-dimensional evolution algebras over the base field $\mathbb{K}$ with basis $\{e_1,\ldots,e_n\}$, whereas $\mathcal{T}_n(\mathbb{K})$ denotes the direct product $\prod_{i=1}^n \langle\,e_1,\ldots,e_n\,\rangle$. Every evolution algebra in $\mathcal{E}_n(\mathbb{K})$ is uniquely determined by an {\em structure tuple} $T=(\mathfrak{t}_1,\ldots,\mathfrak{t}_n)\in \mathcal{T}_n(\mathbb{K})$, where $\mathfrak{t}_i=e_ie_i$, for all $i\leq n$. The structure tuple $T$ also determines the {\em genetic pattern} of the corresponding mitosis process.

\begin{prop}[\cite{Falcon2016_Evol}] \label{prop_0} Let $\mathbb{K}$ be a field. The next results hold.
\begin{enumerate}[a)]
\item Any two structure tuples in $\mathcal{T}_n(\mathbb{K})$ that are equal up to permutation of their components and basis vectors give rise to a pair of strongly isotopic evolution algebras.
\item Let $T$ be a structure tuple in $\mathcal{T}_n(\mathbb{K})$. There always exists a structure tuple $T'=(\sum_{j=1}^nc_{1j}e_j,\ldots,\sum_{j=1}^nc_{nj}e_j)\in \mathcal{T}_n(\mathbb{K})$ such that
\begin{enumerate}[a)]
\item If $c_{ii}=0$, for some $i\geq 1$, then $c_{jk}=0$, for all $j,k\geq i$.
\item If $c_{ii}\neq 0$, for some $i\geq 1$, then $c_{ij}=0$, for all $j\neq i$.
\item The evolution algebra in $\mathcal{E}_n(\mathbb{K})$ of structure tuple $T'$ is strongly isotopic to that one of structure tuple $T$.
\end{enumerate}
\end{enumerate}
\end{prop}

Let $A$ be an evolution algebra in $\mathcal{E}_n(\mathbb{K})$. Isotopisms preserve the dimension of the {\em derived algebra} $A^2=\{uv\mid\, u,v\in A\}\subseteq A$ and that of the {\em annihilator} $\mathrm{Ann}(A)=\{u\in A\mid\, uv=0, \text { for all } v\in A\}$. Hereafter, $\mathcal{E}_{n;m}(\mathbb{K})$ denotes the subset of evolution algebras in $\mathcal{E}_n(\mathbb{K})$ with an $m$-codimensional annihilator.

\begin{theo}[\cite{Falcon2016_Evol}] \label{theo_0} Let $\mathbb{K}$ be a field. The next results hold.
\begin{enumerate}[a)]
\item Let $m<m'$. None evolution algebra in $\mathcal{E}_{n;m}(\mathbb{K})$ is isotopic to an evolution algebra in $\mathcal{E}_{n;m'}(\mathbb{K})$.
\item The set $\mathcal{E}_{n;0}(\mathbb{K})$ is only formed by the $n$-dimensional trivial algebra, which have all its structure constants equal to zero.
\item Any algebra in $\mathcal{E}_{n;1}(\mathbb{K})$ is isotopic to the algebra described as $e_1e_1=e_1$. If the former is not isomorphic to the latter, then it is isomorphic to the evolution algebra described as $e_1e_1=e_2$.
\item Any algebra in $\mathcal{E}_{n;2}(\mathbb{K})$ is isotopic to the algebra described as $e_1e_1=e_2e_2=e_1$ or to that described as $e_1e_1=e_1$ and $e_2e_2=e_2$, and is isomorphic to an evolution algebra in $\mathcal{E}_{n;2}(\mathbb{K})$ such that $e_1e_1\in\{e_1,e_2,e_1+e_2\}$.
\end{enumerate}
\end{theo}

\subsection{Computational Algebraic Geometry}

Let $I$ be an ideal of a multivariate polynomial ring $\mathbb{K}[X]$. The {\em algebraic set} defined by $I$ is the set of common zeros of all its polynomials. The largest monomial of a polynomial in $I$ with respect to a given monomial term ordering is its {\em leading monomial}. The {\em reduced Gr\"obner basis} of $I$ is the only subset $G$ of monic polynomials in $I$ whose leading monomials generate the ideal also generated by all the leading monomials of $I$ and such that no monomial of a polynomial in $G$ is generated by the leading monomials of the rest of polynomials in the basis. This can always be computed from the sequential multivariate division of polynomials described by Buchberger's algorithm \cite{Buchberger2006}, which is extremely sensitive to the number of variables \cite{Lakshman1990a}.

\vspace{0.1cm}

Computational Algebraic Geometry can be used to determine the set of isotopisms and isomorphisms between two evolution algebras $A$ and $A'$ in $\mathcal{E}_n(\mathbb{K})$, with respective structure constants $c_{ij}$ and ${c'}_{ij}$. To this end, we define the sets of variables $\mathfrak{F}_n=\{\mathfrak{f}_{ij}\mid\, i,j\leq n\}$, $\mathfrak{G}_n=\{\mathfrak{g}_{ij}\mid\, i,j\leq n\}$ and $\mathfrak{H}_n=\{\mathfrak{h}_{ij}\mid\, i,j\leq n\}$, which play the role of the entries in the regular matrices related to a possible isotopism $(f,g,h)$ between $A$ and $A'$. The coefficients of each basis vector $e_l$ in the expression $f(e_i)g(e_j)=h(e_ie_j)$ enable us to ensure that $\sum_{k=1}^n \mathfrak{f}_{ik}\mathfrak{g}_{jk}{c'}_{kl}$ is equal to $0$, if $i\neq j$, or to $\sum_{k=1}^n \mathfrak{h}_{kl}c_{ik}$, otherwise.

\begin{theo}[\cite{Falcon2016_Evol}] \label{theo_1} The set of isotopisms and isomorphisms between $A$ and $A'$ are respectively identified with
\begin{enumerate}[a)]
\item The subset of zeros $(f_{11},\ldots,f_{nn},g_{11},\ldots,g_{nn},h_{11},\ldots,h_{nn})\in \mathbb{K}^{3n^2}$ in the algebraic set defined by the ideal $I_{A,A'}$ of $\mathbb{K}[\mathfrak{F}_n\cup\mathfrak{G}_n\cup\mathfrak{H}_n]$ described as
$$\langle\, \sum_{k=1}^n \mathfrak{f}_{ik}\mathfrak{g}_{jk}{c'}_{kl}\mid\, i,j,l\leq n; i\neq j\,\rangle + \langle\, \sum_{k=1}^n \mathfrak{f}_{ik}\mathfrak{g}_{ik}{c'}_{kl} - \sum_{k=1}^n \mathfrak{h}_{kl}c_{ik}\mid\, i,l\leq n\,\rangle,$$
giving rise to non-singular matrices $F=(f_{ij})$, $G=(g_{ij})$ and $H=(h_{ij})$.
\item The subset of zeros $(f_{11},\ldots,f_{nn})\in \mathbb{K}^{n^2}$ in the algebraic set defined by the ideal $J_{A,A'}$ of $\mathbb{K}[\mathfrak{F}_n]$ described as
$$\langle\, \sum_{k=1}^n \mathfrak{f}_{ik}\mathfrak{f}_{jk}{c'}_{kl}\mid\, i,j,l\leq n; i\neq j\,\rangle + \langle\, \sum_{k=1}^n \mathfrak{f}^2_{ik}{c'}_{kl} - \sum_{k=1}^n \mathfrak{f}_{kl}c_{ik}\mid\, i,l\leq n\,\rangle,$$
giving rise to a non-singular matrix $F=(f_{ij})$.
\end{enumerate}
\end{theo}

\vspace{0.5cm}

The computation of the reduced Gr\"obner bases of both ideals in Theorem \ref{theo_1} has recently been implemented \cite{Falcon2016_Evol} in the open computer algebra system for polynomial computations {\sc Singular} \cite{Decker2016}. Particularly, in case of dealing with the finite field $\mathbb{K}=\mathbb{F}_q$, with $q$ a power prime, the respective complexity times that Buchberger's algorithm requires to compute the bases under consideration are $q^{O(3n^2)}+O(n^8)$ and $q^{O(n^2)}+O(n^8)$. We have made use of the mentioned procedure to determine all the isotopisms and isomorphisms that appear throughout this paper. The run time and memory usage that are required to compute each one of the reduced Gr\"obner bases indicated in the paper are respectively $0$ seconds and $0$ MB in a computer system with an {\em Intel Core i7-2600, with a 3.4 GHz processor and 16 GB of RAM}, except for the computation related to Proposition \ref{prop_Isot_0}.d, for which these two measures of computational efficiency have been $3$ seconds and $1$ MB, respectively.

\section{Distribution of $\mathcal{E}_{n;3}(\mathbb{K})$ into isotopism classes}

Let $n\geq 3$ and let $\mathbb{K}$ be a field. In order to determine the distribution of the set $\mathcal{E}_{n;3}(\mathbb{K})$ into isotopism classes, Proposition \ref{prop_0}.b enables us to focus on those $n$-dimensional evolution algebras described as
$$A_{a,b}:=\begin{cases}
e_1e_1=e_1,\\
e_2e_2=\alpha e_1,\\
e_3e_3=\beta e_1,
\end{cases}$$
for some $\alpha,\beta\in\mathbb{K}\setminus\{0\}$, or
$$B_{\alpha,\beta,\gamma}:=\begin{cases}
e_1e_1=e_1,\\
e_2e_2=e_2,\\
e_3e_3=\alpha e_1+\beta e_2+\gamma e_3,
\end{cases}$$
for some $(\alpha,\beta,\gamma)\in\mathbb{K}^3\setminus\{(0,0,0)\}$.

\vspace{0.1cm}

\begin{prop}\label{prop_Isot_0} The next results hold.
\begin{enumerate}[a)]
\item $A_{\alpha,\beta}\sim A_{1,1}$, for all $\alpha,\beta\in\mathbb{K}\setminus\{0\}$.
\item Let $\alpha,\beta,\gamma\in\mathbb{K}$ be such that $\gamma\neq 0$. Then, $B_{\alpha,\beta,\gamma}\sim B_{0,0,1}$.
\item Let $(\alpha,\beta)\in \mathbb{K}^2\setminus\{(0,0)\}$. Then, $B_{\alpha,\beta,0}\sim B_{1,\beta',0}$, for some $\beta'\in\mathbb{K}$.
\item $B_{1,\beta,0}\sim B_{1,1,0}$, for all $\beta\in \mathbb{K}\setminus\{0\}$.
\item The evolution algebras $B_{1,1,0}$ and $B_{1,0,0}$ are not isotopic.
\end{enumerate}
\end{prop}

{\bf Proof.} Let us prove each assertion separately.
\begin{enumerate}[a)]
\item The triple $(f,\mathrm{Id},\mathrm{Id})$ such that $f(e_1)=e_1$, $f(e_2)=\alpha e_2$ and $f(e_3)=\beta e_3$ is an isotopism between the evolution algebras $A_{\alpha,\beta}$ and $A_{1,1}$.
\item Since isotopisms preserve the dimension of derived algebras, Proposition \ref{prop_0}.b establishes that any evolution algebra $B_{\alpha,\beta,\gamma}$ with a three-dimensional derived algebra (that is, such that $\gamma\neq 0$) is isotopic to the evolution algebra $B_{0,0,1}$.
\item From Proposition \ref{prop_0}.a, we can suppose $\alpha\neq 0$. Otherwise, it is enough to switch the vectors $e_1$ and $e_2$. The triple $(f,\mathrm{Id},\mathrm{Id})$ such that $f(e_i)=e_i$, for $i\in\{1,2\}$, and $f(e_3)=\alpha e_3$, is then an isotopism between the evolution algebras $B_{\alpha,\beta,0}$ and $B_{1,\beta/\alpha,0}$.
\item The triple $(f,g,h)$ related to the non-singular matrices
    $$F=\left(\begin{array}{ccc}
    0 & \beta & 0\\
    1 & 0 & 0\\
    0 & 0 & \beta\end{array}\right)\hspace{1cm}
    G=\left(\begin{array}{ccc}
    0 & 1 & 0\\
    1 & 0 & 0\\
    0 & 0 & 1\end{array}\right)\hspace{1cm}
    H=\left(\begin{array}{ccc}
    0 & \beta & 0\\
    1 & 0 & 0\\
    0 & 0 & 1\end{array}\right)$$
is an isotopism between the evolution algebras $B_{1,\beta,0}$ and $B_{1,1,0}$.
\item Let $(f,g,h)$ be an isotopism between the evolution algebras $B_{1,1,0}$ and $B_{1,0,0}$ and let $H=(h_{ij})$ be the non-singular matrix that is related to the linear transformation $h$. The computation of the reduced Gr\"obner basis related to the ideal $I_{B_{1,1,0},B_{1,0,0}}$ in Theorem \ref{theo_1}.a enables us to ensure in particular that $h_{12}=h_{13}=h_{22}=h_{23}=0$. But then, the matrix $H$ is singular, which contradicts the fact of being $(f,g,h)$ an isotopism. Hence, the algebras $B_{1,1,0}$ and $B_{1,0,0}$ are not isotopic. \hfill $\Box$
\end{enumerate}

\begin{theo}\label{theo_Isot_0} There exists four isotopism classes in $\mathcal{E}_{n;3}(\mathbb{K})$, whatever the base field is. They correspond to the evolution algebras $A_{1,1}$, $B_{0,0,1}$, $B_{1,1,0}$ and $B_{1,0,0}$.
\end{theo}

{\bf Proof.} The result follow straightforwardly from Proposition \ref{prop_Isot_0}.\hfill $\Box$

\begin{coro}\label{coro_Isot_0} There exist eight isotopism classes of three-dimensional evolution algebras over any field. They correspond to the evolution algebras of structure tuples $(0,0,0)$, $(e_1,0,0)$, $(e_1,e_1,0)$, $(e_1,e_2,0)$, $(e_1,e_1,e_1)$, $(e_1,e_2,e_3)$, $(e_1,e_2,e_1+e_2)$ and $(e_1,e_2,e_1)$.
\end{coro}

{\bf Proof.} The result is an immediate consequence of Theorems \ref{theo_0} and \ref{theo_Isot_0}. \hfill $\Box$

\vspace{0.1cm}

Let us describe the spectrum of genetic patterns of three distinct genotypes during a mitosis process that are related to the classification exposed in Corollary \ref{coro_Isot_0}. To this end, let $u$, $v$ and $w$ be three distinct elements of the evolution algebra under consideration. The spectrum is formed by the next eight genetic patterns: $(0,0,0)$, for which no offspring exists; $(u,0,0)$, for which only one of the genotypes gives rise to offspring; $(u,u,0)$, for which exactly one of the genotypes does not produce offspring, whereas the other two give rise to offspring with the same genotype; $(u,v,0)$, for which exactly one of the genotypes does not produce offspring, whereas the other two give rise to offspring with distinct genotypes; $(u,u,u)$, for which the offspring has always the same genotype, whatever the initial one is; $(u,v,w)$, for which the genotype of the offspring depends directly on that of the cell parent; $(u,v,\frac 12 u+\frac 12v)$, for which the third genotype gives rise to each one of the genotypes produced by the other two with the same probability; and $(u,v,u)$, for which two of the genotypes produce offspring with the same genotype.

\section{Distribution of three-dimensional evolution algebras with a non-trivial annihilator into isomorphism classes}

Similarly to the reasoning exposed in the previous section, the distribution of the set $\mathcal{E}_3(\mathbb{K})$ into isomorphism classes can be partitioned into that of the subsets $\mathcal{E}_{3;m}(\mathbb{K})$, for $m\leq 3$. In particular, Theorem \ref{theo_0} describes these partitions in case of dealing with $m\in\{0,1\}$. Specifically, the only evolution algebra in $\mathcal{E}_{3;0}(\mathbb{K})$ is the trivial three-dimensional algebra, whereas there exist exactly two isomorphism classes in $\mathcal{E}_{3;1}(\mathbb{K})$. They correspond to the three-dimensional algebra described by $e_1e_1=e_1$ and that one described by $e_1e_1=e_2$. This section deals, therefore, with the computation of isomorphism classes in $\mathcal{E}_{3;2}(\mathbb{K})$. To this end, Theorem \ref{theo_0}.d enables us to focus on those three-dimensional evolution algebras described as
$$C_{\alpha,\beta,\gamma,\delta,\epsilon}:=\begin{cases}
e_1e_1=\alpha e_1+\beta e_2,\\
e_2e_2=\gamma e_1+\delta e_2+\epsilon e_3,
\end{cases}$$
or
$$D_{\alpha,\beta,\gamma,\delta,\epsilon}:=\begin{cases}
e_1e_1=\alpha e_1+\beta e_2,\\
e_3e_3=\gamma e_1+\delta e_2+\epsilon e_3,
\end{cases}$$
for some $(\alpha,\beta)\in\{0,1\}^2\setminus\{(0,0)\}$ and $(\gamma,\delta,\epsilon)\in\mathbb{K}^3\setminus\{(0,0,0)\}$. The next result establishes when an evolution algebra of type $D_{\alpha,\beta,\gamma,\delta,\epsilon}$ is isomorphic to one of type $C_{\alpha,\beta,\gamma,\delta,\epsilon}$. Throughout the distinct proofs of this and the subsequent results, we denote, respectively, as $F=(f_{ij})$ and $G$, the non-singular matrix related to an isomorphism $f$ between the algebras under consideration and the corresponding reduced Gr\"obner basis described in Theorem \ref{theo_1}.b.

\begin{prop}\label{prop_CD} Let $(\alpha,\beta)\in\{0,1\}^2\setminus\{(0,0)\}$ and $(\gamma,\delta,\epsilon)\in\mathbb{K}^3\setminus\{(0,0,0)\}$. Then,
\begin{enumerate}[a)]
\item $D_{1,0,\gamma,\delta,\epsilon}\cong C_{1,0,\gamma,\epsilon,\delta}$
\item $D_{1,1,\gamma,\delta,\epsilon}\cong C_{1,0,\gamma,\epsilon,\delta-\gamma}$.
\item If $\epsilon\neq 0$, then $D_{0,1,\gamma,\delta,\epsilon}\cong  C_{\epsilon,\gamma,0,0,1}$.
\item If $\gamma\neq 0$, then $D_{0,1,\gamma,\delta,0}\cong C_{0,\gamma,0,0,1}$.
\item Let $\delta'\in\mathbb{K}\setminus\{0\}$. Then,
\begin{enumerate}[i.]
\item The evolution algebras $D_{0,1,0,\delta',0}$ and $C_{\alpha,\beta,\gamma,\delta,\epsilon}$ are not isomorphic.
\item If $\delta\neq 0$, then the evolution algebras $D_{0,1,0,\delta,0}$ and $D_{0,1,0,\delta',0}$ are isomorphic if and only if there exists $m\in\mathbb{K}\setminus\{0\}$ such that $\delta=m^2 \delta'$.
\end{enumerate}
\end{enumerate}
\end{prop}

{\bf Proof.} Let us prove each assertion separately.
\begin{enumerate}[a)]
\item It is enough to switch the basis vectors $e_2$ and $e_3$.
\item The result follows from (a) and the fact of being $D_{1,1,\gamma,\delta,\epsilon}\cong D_{1,0,\gamma,\delta-\gamma,\epsilon}$ by means of the isomorphism $f$ such that $f(e_1)=e_1-e_2$ and $f(e_i)=e_i$, for all $i\in\{2,3\}$.
\item It is enough to consider the isomorphism $f$ such that $f(e_1)=e_2$, $f(e_2)=e_3$ and $f(e_3)=e_1-\delta e_3/\epsilon$.
\item Both algebras are isomorphic by means of the isomorphism $f$ such that $f(e_1)=e_2-\delta e_3/\gamma$, $f(e_2)=e_3$ and $f(e_3)=e_1$.
\item \begin{enumerate}[i.]
\item The computation of the related reduced Gr\"obner basis $G$ implies that $\delta'f_{21}^2=\delta'f_{22}^2=\delta'f_{23}^2=0$. Since $\delta'\neq 0$, we have that $f_{21}=f_{22}=f_{23}=0$, which contradicts the fact of being $F$ a non-singular matrix.
\item Here, the corresponding reduced Gr\"obner basis $G$ enables us to ensure that $f_{12}=f_{21}=f_{23}=f_{31}=f_{32}=0$, $f_{22}=f_{11}^2$ and $\delta=f_{22}\delta'$. Hence, the isomorphism $f$ exists if and only if there exists $m\in\mathbb{K}\setminus\{0\}$ such that $\delta=m^2 \delta'$. In such a case, it is enough to consider the linear transformation $f$ such that $f(e_1)=m e_1$, $f(e_2)=m^2 e_2$ and $f(e_3)=e_3$. \hfill $\Box$
\end{enumerate}
\end{enumerate}

\vspace{0.4cm}

Let us focus now on those evolution algebras of type $C_{\alpha,\beta,\gamma,\delta,\epsilon}$. Since isomorphisms preserve the dimension of derived algebras, it is remarkable the possibility of distinguishing whether the parameter $\epsilon$ is equal to $0$ or not. In the first case, the distribution of the set $\mathcal{E}_{3;2}(\mathbb{K})$ into isomorphism classes is equivalent to that of the set $\mathcal{E}_{2;2}(\mathbb{K})$, which has recently been discussed by the authors \cite{Falcon2016_Evol}.

\begin{theo}\label{theo_clev1} Let $(\alpha,\beta)\in\{0,1\}^2\setminus\{(0,0)\}$ and $(\gamma,\delta)\in\mathbb{K}^2\setminus\{(0,0)\}$. The evolution algebra $C_{\alpha,\beta,\gamma,\delta,0}\in\mathcal{E}_{3;2}(\mathbb{K})$ is isomorphic to exactly one of the next algebras
\begin{itemize}
\item $C_{1,0,\gamma,0,0}$, with $\gamma\in\mathbb{K}\setminus\{0\}$. Here, $C_{1,0,\gamma,0,0}\cong C_{1,0,\gamma',0,0}$ if and only if $\gamma' = \gamma m^2$ for some $m\in\mathbb{K}\setminus\{0\}$.
\item $C_{1,1,-1,-1,0}$.
\item $C_{1,0,\gamma,\delta,0}$, with $\delta\neq 0$. Here, $C_{1,0,\gamma,\delta,0}\cong C_{1,0,\gamma',\delta',0}$ if and only if $\gamma \delta'^2 = \delta^2\gamma$.
\item $C_{0,1,\gamma,\delta,0}$, with $\gamma\neq 0$. Here, $C_{0,1,\gamma,\delta,0}\cong C_{0,1,\gamma',\delta',0}$ if and only if there exists an element $m\in\mathbb{K}\setminus\{0\}$ such that $\gamma = \gamma' m^3$ and $\delta=\delta' m^2$, or, $\gamma=\gamma'^2 m^3$ and $\delta=\delta'=0$.
\item $C_{1,1,\gamma,\delta,0}$, with $\gamma\neq d$. Here, if $\gamma\neq 0\neq \delta$, then $C_{1,1,\gamma,\delta,0}\cong C_{1,1,\gamma',\delta',0}$ if and only if $\gamma'=\gamma^2/\delta^3$ and $\delta'=\gamma/\delta^2$.
\end{itemize}
\end{theo}

\vspace{0.1cm}

We focus here, therefore, on the second case.

\begin{lemm}\label{lemm_Isom_0} Let $(\alpha,\beta)\in\{0,1\}^2\setminus\{(0,0)\}$ and let $\gamma,\delta,\epsilon\in\mathbb{K}$ be such that $\epsilon\neq 0$. Then, $C_{\alpha,\beta,\gamma,\delta,\epsilon}\cong C_{\alpha,\beta,\gamma,\delta,1}$.
\end{lemm}

{\bf Proof.} It is enough to consider the linear transformation $f$ between both algebras that is described as $f(e_i)=e_i$, if $i\in\{1,2\}$, and $f(e_3)=\frac 1\epsilon e_3$. \hfill $\Box$

\vspace{0.1cm}

Lemma \ref{lemm_Isom_0} enables us to focus on the study of evolution algebras of the form $C_{\alpha,\beta,\gamma,\delta,1}$. The next result establishes which ones of these algebras are isomorphic to the evolution algebras described in Theorem \ref{theo_clev1}.

\begin{prop}\label{prop_Isom_0} Let $(\alpha,\beta)\in\{0,1\}^2\setminus\{(0,0)\}$ and let $\gamma,\delta\in\mathbb{K}$. Then,
\begin{enumerate}[a)]
\item The evolution algebra $C_{\alpha,\beta,\gamma,\delta,1}$ is not isomorphic to any algebra $C_{1,0,\gamma',0,0}$, with $\gamma'\in\mathbb{K}\setminus\{0\}$, or to $C_{1,1,-1,-1,0}$.
\item Let $\gamma',\delta'\in\mathbb{K}$ be such that $\delta'\neq 0$. The evolution algebra $C_{\alpha,\beta,\gamma,\delta,1}$ is isomorphic to the algebra $C_{1,0,\gamma',\delta',0}$ only if $\beta=0$. Specifically,
    \begin{enumerate}[i.]
    \item $C_{1,0,\gamma,\delta,1}\cong C_{1,0,\gamma,\delta,0}$, whenever $\delta\neq 0$.
    \item The evolution algebra $C_{1,0,\gamma,0,1}$ and $C_{1,0,\gamma',\delta',0}$ are not isomorphic.
    \end{enumerate}
\item Let $\gamma',\delta'\in\mathbb{K}$ be such that $\gamma'\neq 0$. Then,
    \begin{enumerate}[i.]
    \item The evolution algebras $C_{1,0,\gamma,0,1}$ and $C_{0,1,\gamma',\delta',0}$ are not isomorphic.
    \item $C_{0,1,\gamma,\delta,1}\cong C_{0,1,\gamma,\delta,0}$, whenever $\gamma\neq 0$.
    \item The evolution algebra $C_{0,1,0,\delta,1}$ and $C_{0,1,\gamma',\delta',0}$ are not isomorphic.
    \item The evolution algebra $C_{1,1,\gamma,\delta,1}$ is isomorphic to the algebra $C_{0,1,\gamma',\delta',0}$ if and only if $\delta=0$ and $\gamma'^2=\gamma\delta'^3$. In particular, $C_{1,1,\gamma,0,1}\cong C_{0,1,\gamma^2,\gamma,0}$, for all $\gamma\in\mathbb{K}\setminus\{0\}$.
    \end{enumerate}
\item Let $\gamma',\delta'\in\mathbb{K}$ be such that $\gamma'\neq \delta'$. Then,
\begin{enumerate}[i.]
    \item None of the evolution algebras $C_{1,0,\gamma,0,1}$, $C_{0,1,0,\delta,1}$ and $C_{1,1,0,0,1}$ is isomorphic to the algebra $C_{1,1,\gamma',\delta',0}$.
    \item Let $\delta\neq 0$. Then,
     \begin{itemize}
     \item $C_{1,1,\gamma,\delta,1}\cong C_{1,1,\gamma^2/\delta^3,\gamma/\delta^2,0}$ if and only if $\delta\neq \gamma\neq 0$.
     \item $C_{1,1,0,\delta,1}\cong C_{1,1,0,\delta,0}$.
     \end{itemize}
\end{enumerate}
\end{enumerate}
\end{prop}

{\bf Proof.} Let us prove each assertion separately.
\begin{enumerate}[a)]
\item The result follows straightforwardly from the fact that isomorphisms preserve the dimension of derived algebras. In the case under consideration, the derived algebra of $C_{\alpha,\beta,\gamma,\delta,1}$ is two-dimensional, whereas those of $C_{1,0,\gamma',0,0}$ and $C_{1,1,-1,-1,0}$ are one-dimensional.
\item The computation of the corresponding reduced Gr\"obner basis $G$ enables us to ensure that $f_{31}\delta'=f_{32}\delta'=0$ and hence, $f_{31}=f_{32}=0$, because $\delta'\neq 0$. After imposing this condition to the polynomials in $G$, we obtain that $f_{11}f_{22}=f_{12}f_{21}=f_{21}\beta\gamma=f_{11}\beta\gamma=0$. From the non-singularity of the matrix $F$, we have, therefore, that $\beta\gamma=0$. Now, if $\beta=1$ and $\gamma=0$, then the basis $G$ involves that $f_{11}=f_{22}=0$ and $f_{21}f_{33}=0$, which is a contradiction with being $F$ a non-singular matrix. Hence, we can focus on the case $\beta=0$.
    \begin{enumerate}[i.]
    \item If $\delta\neq 0$, then the evolution algebras $C_{1,0,\gamma,\delta,1}$ and $C_{1,0,\gamma,\delta,0}$ are isomorphic by means of the isomorphism related to the non-singular matrix
    $$\left(\begin{array}{ccc}
    1 & 0 & 0\\
    0 & 1 & 1\\
    0 & 0 & -\delta
    \end{array}\right).$$
    \item The computation of the basis $G$ involves that $f_{22}\delta'=0$ and hence, $f_{22}=0$, because $\delta'\neq 0$. After imposing this condition to the polynomials in $G$, we also obtain that $f_{11}=f_{21}=0$, which is a contradiction with being $F$ a non-singular matrix.
    \end{enumerate}
\item \begin{enumerate}[i.]
    \item The computation of the basis $G$ involves in this case that $f(e_3)=0$, which is a contradiction with being $F$ a non-singular matrix.
    \item If $\gamma\neq 0$, then the evolution algebras $C_{0,1,\gamma,\delta,1}$ and $C_{0,1,\gamma,\delta,0}$ are isomorphic by means of the isomorphism related to the non-singular matrix
    $$\left(\begin{array}{ccc}
    1 & 0 & 1\\
    0 & 1 & 0\\
    0 & 0 & -\gamma
    \end{array}\right).$$
    \item The computation of the basis $G$ enables us to ensure that $f_{11}=f_{12}=f_{31}=f_{32}=0$, which contradicts the fact of being $F$ a non-singular matrix.
    \item The computation of the basis $G$ involves now that $f_{11}=f_{22}=f_{31}=f_{32}=f_{12}\delta=0$. From the non-singularity of the matrix $F$, we have that $\delta=0$. After imposing these conditions to the polynomials in $G$, we obtain that the isotopism $f$ exists if and only $\gamma'^2=\gamma \delta'^3$. In such a case, $\gamma\neq 0\neq \delta'$ and it is enough to define $f$ as the linear transformation related to the non-singular matrix
    $$\left(\begin{array}{ccc}
    0 & 1/\delta' & 1\\
    \gamma'/\delta'^2 & 0 & -1\\
    0 & 0 & -\gamma
    \end{array}\right).$$
    \end{enumerate}
\item
    \begin{enumerate}[i.]
    \item In the three cases, the computation of the corresponding reduced Gr\"obner basis $B$ involves that $f(e_3)=0$, which contradicts the fact of being $F$ a non-singular matrix.
    \item The computation of the basis $G$ enables us to ensure in this case that $f_{31}(\gamma'-\delta')=0$ and hence, $f_{31}=0$, because $\gamma'\neq \delta'$. After imposing this condition to the polynomials in $G$, we obtain that $f_{32}\delta'=0$ and hence, $f_{32}=0$, because $\delta'\neq 0$. The non-singularity of the matrix $F$ involves then that $f_{33}=f_{13}(\delta-\gamma)\neq 0$, with $\delta\neq \gamma$. Then, from the generators in $G$, exactly one of the next two cases holds.
        \begin{itemize}
        \item $f_{11}=f_{22}=0$, $\gamma=f_{21}\neq 0$ and $f_{12}=1/\delta'$. According to the basis $G$, the isomorphism $f$ exists in this case if and only if $\gamma'=\gamma^2/\delta^3$ and $\delta'=\gamma/\delta^2$. In particular, it must be $\gamma\neq 0$. In such a case, it is enough to define $f$ as the linear transformation related to the non-singular matrix
    $$\left(\begin{array}{ccc}
    0 & \delta^2/\gamma & 1\\
    \delta  & 0 & -1\\
    0 & 0 & \delta-\gamma
    \end{array}\right).$$
        \item $f_{12}=f_{21}=0$ and $f_{11}=f_{22}=1$. In this case, the linear transformation $f$ related to the non-singular matrix
    $$\left(\begin{array}{ccc}
    1 & 0 & 1\\
    0  & 1 & -1\\
    0 & 0 & \delta
    \end{array}\right)$$
    is an isomorphism between $C_{1,1,0,\delta,1}$ and $C_{1,1,0,\delta,0}$. \hfill $\Box$
        \end{itemize}
    \end{enumerate}
\end{enumerate}

\vspace{0.5cm}

Let us focus now on those evolution algebras of the form $C_{\alpha,\beta,\gamma,\delta,1}$ that are not isomorphic to any of the algebras described in Theorem \ref{theo_clev1}, that is, on the evolution algebras of the form $C_{1,0,\gamma,0,1}$, $C_{0,1,0,\delta,1}$ and $C_{1,1,\gamma,\gamma,1}$.

\begin{lemm}\label{lem_Isom_1_0} Let $(\alpha,\beta)$ and $(\alpha',\beta')$ be two pairs in $\{0,1\}^2\setminus\{(0,0)\}$ and let $\gamma,\delta,\gamma',\delta'\in\mathbb{K}$ be such that the evolution algebras $C_{\alpha,\beta,\gamma,\delta,1}$ and $C_{\alpha',\beta',\gamma',\delta',1}$ are isomorphic. Then, $f(e_3)=m e_3$, for some $m\in \mathbb{K}\setminus\{0\}$.
\end{lemm}

{\bf Proof.} The coefficients of the basis vector $e_3$ in the expression $f(e_3)f(e_3)$ $=f(e_3e_3)=0$ involves that $f_{32}=0$. Due to it, the coefficients of the basis vectors $e_1$ and $e_2$ in the previous expression enable us to ensure that $f_{31}\alpha'=f_{31}\beta'=0$. Since $(\alpha',\beta')\neq (0,0)$, it is $f_{31}=0$. Hence, $f(e_3)=f_{33}e_3$, where $f_{33}\neq 0$, because the matrix $F$ is non-singular.\hfill $\Box$

\vspace{0.1cm}

\begin{prop}\label{prop_Isom_1} Let $\gamma,\delta,\gamma'\in\mathbb{K}$. Then,
\begin{enumerate}[a)]
\item The evolution algebras $C_{1,0,\gamma,0,1}$ and $C_{1,0,\gamma',0,1}$ are isomorphic if and only if there exists $m\in\mathbb{K}\setminus\{0\}$ such that $\gamma=m^2 \gamma'$.
\item $C_{0,1,0,\delta,1}\cong C_{1,0,\delta,0,1}$, whenever $\delta\neq 0$.
\item The evolution algebras $C_{0,1,0,0,1}$ and $C_{1,1,\gamma,\gamma,1}$ are not isomorphic. Besides, none of them is isomorphic to the algebra $C_{1,0,\gamma',0,1}$.
\item If $\gamma\neq\gamma'$, then the evolution algebras $C_{1,1,\gamma,\gamma,1}$ and $C_{1,1,\gamma',\gamma',1}$ are isomorphic if and only if $\gamma\gamma'=1$.
\end{enumerate}
\end{prop}

{\bf Proof.} Let us prove each assertion separately.
\begin{enumerate}[a)]
\item The computation of the corresponding reduced Gr\"obner basis $G$ enables us to ensure that $f_{11}=1$, $f_{12}=f_{13}=f_{21}=f_{31}=f_{32}=0$, $f_{33}=f_{22}^2$ and $\gamma=f_{33}\gamma'$. From Lemma \ref{lem_Isom_1_0}, the isomorphism $f$ exists if and only if there exists $m\in\mathbb{K}\setminus\{0\}$ such that $\gamma=m^2 \gamma'$. In such a case, it is enough to consider the linear transformation $f$ such that $f(e_1)=e_1$, $f(e_2)=m e_2$ and $f(e_3)=m^2 e_3$.
\item The linear transformation related to the non-singular matrix $$\left(\begin{array}{ccc}
    0 & 1 & 0\\
    \delta & 0 & 1\\
    0 & 0 & -\delta
    \end{array}\right)$$
    is an isomorphism between these two algebras.
\item In any of the possible cases under consideration, the computation of the corresponding reduced Gr\"obner basis $G$ involves that $f_{33}=0$, which contradicts Lemma \ref{lem_Isom_1_0}.
\item The computation of the corresponding reduced Gr\"obner basis $G$ involves that $f_{11}=f_{22}$, $f_{12}=f_{21}=\gamma (1-f_{22})$ and $f_{31}=f_{32}=0=(\gamma'-\gamma) f_{22}$. Hence, since $\gamma\neq \gamma'$, we have that $f_{11}=f_{22}=0$ and $f_{12}=f_{21}=\gamma$. After imposing these conditions to the polynomials in $G$, we obtain that $\gamma' f_{33}=-\gamma^2$ and $\gamma^2\gamma'=\gamma$. The first condition, together with Lemma \ref{lem_Isom_1_0}, enables us to ensure that $\gamma\neq 0\neq \gamma'$. Then, the second condition involves that $\gamma\gamma'=1$ and hence, $f_{33}=-\gamma^3$. We also obtain from the reduced Gr\"obner basis $G$ that $f_{13}+f_{23}=\gamma^2$. In particular, the linear transformation related to the non-singular matrix $$\left(\begin{array}{ccc}
    0 & \gamma & \gamma^2\\
    \gamma & 0 & 0\\
    0 & 0 & -\gamma^3
    \end{array}\right)$$
    is an isomorphism between the evolution algebras $C_{1,1,\gamma,\gamma,1}$ and $C_{1,1,\gamma',\gamma',1}$, whenever $\gamma\gamma'=1$.\hfill $\Box$
\end{enumerate}

\vspace{0.5cm}

The next theorem gathers together all the results that we have exposed throughout this section.

\begin{theo}\label{theo_Isom_0} Any evolution algebra in the set $\mathcal{E}_{3;2}(\mathbb{K})$ is isomorphic to exactly one of the algebras of Theorem \ref{theo_clev1} or to one of the next algebras
\begin{itemize}
\item $C_{1,0,\gamma,0,1}$, for some $\gamma\in\mathbb{K}$. Here, $C_{1,0,\gamma,0,1}\cong C_{1,0,\gamma',0,1}$ if and only if there exists $m\in\mathbb{K}\setminus\{0\}$ such that $\gamma=m^2 \gamma'$.
\item $C_{1,1,\gamma,\gamma,1}$, for some $\gamma\in\mathbb{K}$. Here, $C_{1,1,\gamma,\gamma,1}\cong C_{1,1,\gamma',\gamma',1}$ if and only if $\gamma\gamma'=1$.
\item $D_{0,1,0,\delta,0}$, for some $\delta\in\mathbb{K}\setminus\{0\}$. Here, $D_{0,1,0,\delta,0}\cong D_{0,1,0,\delta',0}$ if and only if there exists $m\in\mathbb{K}\setminus\{0\}$ such that $\delta=m^2 \delta'$.
\end{itemize}
\end{theo}

\section{Conclusion and further studies}

This paper has dealt with the use of Computational Algebraic Geometry to determine the distribution of three-dimensional evolution algebras over any field into isotopism and isomorphism classes. We have proved in particular that the set $\mathcal{E}_3(\mathbb{K})$ is distributed into eight isotopism classes, whatever the base field is, and we have characterized their isomorphism classes in case of dealing with algebras with a one-dimensional annihilator. We have also described the spectrum of genetic patterns of three distinct genotypes during mitosis of eukaryotic cells. The distribution into isotopism classes of the set $\mathcal{E}_4(\mathbb{K})$ and the characterization of their isomorphism classes in case of dealing with algebras with one-dimensional annihilator is established as further work.

\section*{Acknowledgements}

\vspace{-0.25cm}

The authors are very grateful to the Editors and to the referees for their valuable comments and suggestions that have improved the original version of the paper. Particularly, the authors gratefully thank to one of the referees for the suggestion of considering for our study the recent work developed by Cabrera et al. \cite{Cabrera2017}, which has enabled us to detect a pair of omitted cases of study in the preliminary version of the paper.

\end{document}